\documentclass[12pt]{article}
\def\P{\mathbb{P}}
\def\E{\mathbb{E}}
\def\F{{\cal F}}

\newcommand{\R}        {{{\rm I \hskip -2pt R}}}

\newtheorem{lemma}{Lemma}

\newtheorem{theorem}{Theorem}

\newtheorem{Corollary}{Corollary}

\newtheorem{Remark}{Remark}

\usepackage{amsmath}
\usepackage{amssymb}
\usepackage{color}

\begin{document}
\title{{\normalsize\tt\hfill\jobname.tex}\\
Positive recurrence of a solution of an SDE 
with 
variable switching intensities
}
\author{Alexander Veretennikov\footnote{
Institute for Information Transmission Problems, Moscow, Russian Federation, email: {ayv@iitp.ru}}\;\footnote{This research was funded by Russian Foundation for Basic Research grant 20-01-00575a.
}
}

\maketitle       

\begin{abstract}
Positive recurrence of a $d$-dimensional diffusion 
with switching  and with one recurrent and one transient regimes and variable switching intensities is established under suitable conditions. The approach is based on embedded Markov chains. 

~

{\em Keywords: {diffusion, switching, variable switching intensities, positive recurrence}}

~

{\em MSC codes: 60H10, 60J60}

\end{abstract}

\section{Introduction}
Let us consider the process $(X_t,Z_t)$ with a continuous component $X$ and discrete one $Z$ described by the stochastic differential equation  in $\R^d$
\begin{align}\label{sde}
dX_{t} =b(X_{t}, Z_t)\, dt+ dW_{t}, \quad t\ge 0, 
\quad X_{0} =x, \; Z_0=z, 
\end{align}
for the component $X$, while $Z_t$ is a continuous-time conditionally Markov process given $X$ on the state space $S= \{0,1\}$ with positive intensities of respective transitions  $\lambda_{01}(x) =: \lambda_0(x), \, \& \, \lambda_{10}(x) =: \lambda_1(x)$; here the variable $x$ signifies a certain (arbitrary Borel measurable) dependence on the component $X$; the trajectories of $Z$ are assumed to be c\`adl\`ag; the probabilities of jumps for $Z$ are conditionally independent given the trajectory of the component $X$ (see the precise description in what follows). Denote 
$$
b(x,0) = b_-(x), \quad b(x,1) = b_+(x), 
$$
$$
\overline{\lambda}_0:= \sup_{x,z} \lambda_0(x), 
\quad 
\underline{\lambda}_0:= \inf_{x,z} \lambda_0(x), \quad 
 \overline{\lambda}_1:= \sup_{x,z} \lambda_1(x), 
\quad \underline{\lambda}_1 := \inf_{x,z} \lambda_1(x).
$$
It is assumed that 
\begin{align}\label{lambdaass}
0< \underline{\lambda}_0\wedge
\underline{\lambda}_1 \le \overline{\lambda}_0 \vee \overline{\lambda}_1 < \infty.
\end{align}
These conditions along with the boundedness of the function $b$ in $x$ suffice for the process $(X_t,Z_t)$ to be well-defined.
A rigorous construction of the system $(X,Z)$ of this type may be given by the SDE system
\begin{align}\label{sdexz}
dX_{t} =b(X_{t}, Z_t)\, dt+ 
dW_{t}, \; t\ge 0, \; X_{0} =x \in \mathbb R^d,
 \nonumber \\  \\ \nonumber
dZ_t=1(Z_t=0)d\pi^0_t - 1(Z_t=1)d\pi^1_t,  \; \; Z_0 \in \{0,1\},
\end{align}
where $\pi^i_t, \, i=0,1,$ are two Poisson processes with intensities $\lambda_i(X_t), \, i=0,1$, respectively. More precisely,
$$
\pi^i_t = \bar\pi^i_{\phi_i(t)}, 
$$
where $\bar\pi^i_t, i=0,1$, are, in turn, two standard Poisson processes with a constant intensity one, independent of the Wiener process $(W_t)$ and of each other, and  the  time changes 
$$
t\mapsto \phi_i(t):=\int_0^t \lambda_i(X_s)ds, \quad i=0,1,
$$ 
are applied to each of them, respectively. 

By virtue of the assumption (\ref{lambdaass}) 
the equation between the jumps only concerns the diffusion part of the SDE (\ref{sdexz}), for which it is well-known since \cite{Ver79} that the equation has a pathwise unique strong solution. The jump moments are stopping times with respect to the filtration $(\F_t = \F_t^{W,\pi^0,\pi^1}, \, t\ge 0)$, 
and the position of the system after any jump $(X_\tau,Z_\tau)$ is uniquely determined by the left limiting values $(X_{\tau-},Z_{\tau-})$: 
$$
X_{\tau} = X_{\tau-}, \quad 
Z_{\tau} = 1(Z_{\tau-}=0).
$$
After any such jump, the diffusion part of the SDE is solved starting from the position $X_{\tau}$ until the next jump, say, $\tau'$, of the component $Z$, and the moment of this next jump is determined by the trajectories of $\pi^0_t$ and (or) of $\pi^1_t$ and by the intensity $\lambda_{Z_s}(X_s), s < \tau'$.
Since there might be only a finite numbers of jumps on any bounded  interval of time, then pathwise (and, hence, also weak) uniqueness follows on $[0,\infty)$. Therefore, the process $(X,Z)$ exists and is markovian. Its (quasi-) generator has a form
$$
Lh(x,z) = \frac12\Delta_{x} h(x,z) + b_{z}(x) \nabla_{x} h(x,z) + \lambda_{z}(x) \, (h(x,\bar z) - h(x,z)),
$$
where $\bar z : = 1(z=0)$ (that is, $\bar z$ is not a $z$, the other state from $\{0,1\}$).

~

For any $t>0$ fixed let us define the function
$$
v(s,x,z):= \E_{s,x,z} f(X_t,Z_t).
$$
The vector-function $v(s,x) = (v(s,x,0),v(s,x,1))$ satisfies the system of PDEs
\begin{align*}
v_s(s,x,0) + L^0 v(s,x,0) + \lambda_{0}(x) \, (v(s,x,1) - v(s,x,0))= 0, \quad v(t,x,0)=f(x,0), 
 \\\\
v_s(s,x,1) + L^1 v(s,x,1) + \lambda_{1}(x) \, (v(s,x,0) - v(s,x,1))= 0, \quad v(t,x,1)=f(x,1), 
\end{align*}
where
$$
L^i = \frac12 \,\Delta_x + \langle b(x,i), \nabla_{x}\rangle, \quad i=0,1.
$$
Due to the results of \cite[Theorem 5.5]{Solo} (see also  \cite{KrylovSafonov}), its solution is continuous 
in the variable $s$ for any bounded and continuous $f$. Hence, the process is Feller's (that is, $\E_{x,z} h(X_t)$ is continuous in $x$ and, of course, bounded for any $h\in C_b$ and any $t>0$). Since the process is Markov and c\`adl\`ag, then it is also strong Markov according to the well-known sufficient condition.

~

The SDE solution is assumed ergodic under the regime $Z=0$ and transient under $Z=1$. We are looking for sufficient conditions for positive recurrence of the strong Markov process $(X_t,Z_t)$. 
Such a problem was considered in \cite{Hairer} for the exponentially recurrent case; for other references see 
\cite{Anulova}, \cite{Khasminskii12}, \cite{Mao}, \cite{ShaoYuan}, and the references therein. Under weak ergodic and transient conditions the setting was earlier investigated in \cite{Ver21} for the case of the constant intensities $\lambda_{0}, \lambda_{1}$ (i.e., not depending on $x$). In \cite{PinskyScheutzow} and \cite{PinskyPinsky} other interesting results about the transience and recurrence for diffusions with switching were established; in particular, examples were given of conditions under which the solutions of the SDEs {\em on the half line with reflection} with two transient ``pure'' regimes are recurrent, and, vice versa, where the solutions of the SDEs  with two recurrent ``pure'' regimes are transient. Here we tackle the general case of a combination of one transient and one recurrent regime. In the case of \cite{Ver21} the lengths of intervals between successive jumps of the discrete component were all independent and independent of the Wiener process. In the general case under the consideration in the present paper they are not independent of $W$ via the component $X$, and, hence, not independent of each other. This difficulty will be overcome with the help of certain comparison arguments.

%

~

The paper consists of the sections: Introduction, Main result, Auxiliaries, and Proof of main result.

\section{Main result}

\begin{theorem}\label{thm1}
Let 
the drift $b = (b_+,b_-)$ be bounded and Borel measurable, and let there exist $r_-, r_+,M>0$ such that 
\begin{equation}\label{al}
0< \underline{\lambda}_0 \wedge \underline{\lambda}_1 
\le \bar \lambda_0 \vee \bar \lambda_1 < \infty,
\end{equation}
\begin{equation}\label{b}
x b_-(x) \le - r_-, \quad x b_+(x)\le + r_+, \quad \forall \,|x|\ge M, 
\end{equation}
and 
\begin{equation}\label{c1}
2r_- > d \; 
\quad \& \quad \underline{\lambda}_1(2r_- -d) > \overline{\lambda}_0(2r_++d).
\end{equation}
Then the process $(X,Z)$ is positive recurrent; moreover, 
there exists $C>0$ such that for all $M_1$ large enough and all $x \in \mathbb R^d$ and for $z=0,1$
\begin{equation}\label{e3}
\mathbb E_{x,z}\tau_{M_1} 
\le C (x^2 + 1), 
\end{equation} 
where 
$$
\tau_{M_1} := \inf(t\ge 0:\, |X_t|\le M_1). 
$$
Moreover, the process $(X_t,Z_t)$ has a unique invariant measure, and for each nonrandom initial condition $x,z$ there is a convergence to this measure in total variation when $t\to\infty$.
\end{theorem}

\section{Auxiliaries}
\label{S:2}
Denote $\|b\| = \sup_{x,z}|b(x,z)|$.
Let $M_1 \gg M$ (the value $M_1$ will be specified later). 
Let
$$
T_0 : = \inf(t\ge 0: Z_t = 0), 
$$
and   
$$0 \le T_0 < T_1 < T_2 < \ldots, $$
where $T_n$ for each $n\ge 1$ is defined by induction as
 
$$
T_{n}:= \inf(t>T_{n-1}: \, Z_{T_{n}}-Z_{T_{n}-} \neq 0). 
$$
Let 
$$
\tau : = \inf(T_n\ge 0: \, |X_{T_n}|\le M_1). 
$$
To prove the theorem it suffices to evaluate from above the value $\mathbb E_{x,z}\tau$ because $\tau_{M_1} \le \tau$. Let $\epsilon>0, q<1$ be positive values satisfying the equality 
\begin{equation}\label{lle}
\overline\lambda_0(2r_++d +\epsilon) = q \underline\lambda_1(2r_- -d - \epsilon)
\end{equation}
(see (\ref{c1})). In the proof of the theorem it suffices to assume $|x| > M$. 

\begin{lemma}\label{lem1}
Under the assumptions of the theorem for any $\delta >0$ there exists $M_1$ such that 
\begin{eqnarray}\label{eps}
\max\left[\sup_{|x|>M_1}\mathbb E_{x,z} \left(\!\!\int_0^{T_1}1(\inf_{0\le s\le t}|X_s| \!\le\! M)dt|Z_0\!=\!0\!\right)\!, 
 \right. \nonumber\\ \\ \nonumber\left.
 \!  \sup_{|x|>M_1}\mathbb E_{x,z} \!\left(\!\int_0^{T_0}\!1(\inf_{0\le s\le t}|X_s| \!\le\! M)dt|Z_0\!=\!1\!\right)\right]
<\delta .
\end{eqnarray}
\end{lemma}
Let us denote by $X^i_t, \, i=0,1$ the solutions of the equations 
\begin{equation}\label{sde0}
dX^i_{t} =b(X^i_{t}, i)\, dt+ 
dW_{t}, \quad t\ge 0,
\quad X^i_{0} =x. 
\end{equation}
{\em Proof of lemma 1.} 
Let $Z_0=0$; then $T_0=0$. We have,
$$
\P_{x,0}(X_t = X^0_t, \, 0\le t\le T_1) = 1, 
$$
due to the uniqueness of solutions of the SDEs (\ref{sde}) (or (\ref{sdexz})) and (\ref{sde0}) and because of the property of stochastic integrals \cite[Theorem 2.8.2]{Kry} to coincide almost surely (a.s.) on the set where the integrands are equal.
Therefore, 
we estimate for any $|x| > M$ with $z=0$:
\begin{align*}
\mathbb E_{x,z} \left(\int_0^{T_1}1(\inf_{0\le s\le t}|X_s| \le M)dt|Z_0=0\right)  = 
\mathbb E_{x,z} \int_0^{T_1}1(\inf_{0\le s\le t}|X^0_s| \le M)dt
 \\\\
= \mathbb E_{x,z} \int_0^\infty 1(t<{T_1})1(\inf_{0\le s\le t}|X^0_s| \le M)dt
= \int_0^\infty \mathbb E_{x,z} 1(t<{T_1})1(\inf_{0\le s\le t}|X^0_s| \le M)dt
 \\\\
\stackrel{\forall \, t_0>0}= \int_0^{t_0} \mathbb E_{x,z} 1(t<{T_1})1(\inf_{0\le s\le t}|X^0_s| \le M)dt 
+ \int_{t_0}^\infty \mathbb E_{x,z} 1(t<{T_1}) 
1(\inf_{0\le s\le t}|X^0_s| \le M)dt
 \\\\
\le \int_0^{t_0} \mathbb E_{x,z} 1(\inf_{0\le s\le t}|X^0_s| \le M)dt 
+ \int_{t_0}^\infty \mathbb E_{x,z} 1(t<{T_1}) 
dt
 \\\\
\le t_0 \mathbb P_{x,z} (\inf_{0\le s\le t_0}|X^0_s| \le M) 
+ \int_{t_0}^\infty \exp(-\underline{\lambda}_0 t)
dt.
\end{align*}
Let us fix some $t_0$, so that
$$
t_0 >  -\underline{\lambda}_0^{-1} \, \ln(\underline{\lambda}_0^{}\delta/2).
$$
Then 
$$
\int_{t_0}^\infty  e^{-\underline{\lambda}_0 s}ds <\delta /2.
$$
Now, with this $t_0$ already fixed, by virtue of the boundedness of $b$ there exists $M_1>M$ such that  for any $|x|\ge M_1$ we get
$$
t_0 \, \mathbb P_{x,z}(\inf_{0\le s\le t_0}|X^0_s| \le M) <\delta /2.
$$
Similarly, the bound for the second term in (\ref{eps}) follows if we replace the process $X^0$ by $X^1$ and the intensity $\lambda_0$ by $\lambda_1$.
\hfill 
{\em QED}

\begin{lemma}\label{lem2}
If $M_1$ is large enough, then under the assumptions of the theorem 
for any $|x|>M_1$ for any $k=0,1,\ldots$
\begin{eqnarray}
\mathbb E_{x,z} (X_{T_{2k+1}\wedge   \tau}^2|Z_0=0, {\mathcal F}_{T_{2k}}) 
\le \mathbb E_{x,z} (X_{T_{2k}\wedge   \tau}^2|Z_0=0, {\mathcal F}_{T_{2k}}) 
 \nonumber \\ \label{ele2a} \\
- 1(\tau > T_{2k}) \E(T_{2k+1}\wedge   \tau - T_{2k}\wedge   \tau |Z_0=0, {\mathcal F}_{T_{2k}})((2r_--d)- \epsilon)
 \nonumber \\ \nonumber \\ \label{ele2a0}
\le \mathbb E_{x,z} (X_{T_{2k}\wedge   \tau}^2|Z_0=0, {\mathcal F}_{T_{2k}})- 1(\tau > T_{2k})\bar\lambda_0^{-1}((2r_--d)- \epsilon), 
 \\ \nonumber
 \\ \label{ele2b}
\mathbb E_{x,z} (X_{T_{2k+2}\wedge   \tau}^2|Z_0=1, {\mathcal F}_{T_{2k+1}})  
\le  \mathbb E_{x,z} (X_{T_{2k+1}\wedge   \tau}^2  |Z_0=1, {\mathcal F}_{T_{2k+1}}) 
 \nonumber\\ \\\label{ele2b0}
+  1(\tau > T_{2k+1}) \E(T_{2k+2}\wedge   \tau-T_{2k+1}\wedge   \tau |Z_0=1, {\mathcal F}_{T_{2k+1}})((2r_-+d)+ \epsilon)
 \nonumber \\ \nonumber \\ \label{ele2B0}
\le \mathbb E_{x,z} (X_{T_{2k+1}\wedge   \tau}^2  |Z_0=1, {\mathcal F}_{T_{2k+1}}) 
+ 1(\tau > T_{2k+1}) \underline{\lambda}_1^{-1}((2r_++d)+ \epsilon).
\end{eqnarray}
\end{lemma}

\begin{Corollary}\label{Cor2}
If $M_1$ is large enough, then under the assumptions of the theorem 
for any $|x|>M_1$ for any $k=0,1,\ldots$
\begin{eqnarray*}
\mathbb E_{x,0} X_{T_{2k+1}\wedge   \tau}^2 - 
\mathbb E_{x,0}  X_{T_{2k}\wedge   \tau}^2
 \\\\
\le - \mathbb E_{x,0} 1(\tau > T_{2k}) \E_{x,0}(T_{2k+1}\wedge   \tau - T_{2k}\wedge   \tau |{\mathcal F}_{T_{2k}})((2r_--d)- \epsilon)
 \\\\ 
= - \mathbb E_{x,0} (T_{2k+1}\wedge   \tau - T_{2k}\wedge   \tau)((2r_--d)- \epsilon)
 \\\\
\le - \mathbb E_{x,0} 1(\tau > T_{2k})\bar\lambda_0^{-1}((2r_--d)- \epsilon),
\end{eqnarray*}
and
\begin{eqnarray*}
\mathbb E_{x,1} X_{T_{2k+2}\wedge   \tau}^2 - \mathbb E_{x,1}  X_{T_{2k+1}\wedge   \tau}^2  
 \nonumber\\ \\
\le  \E_{x,1} 1(\tau > T_{2k+1}) (T_{2k+2}\wedge   \tau-T_{2k+1}\wedge   \tau ) ((2r_-+d)+ \epsilon)
 \nonumber\\ \\
=  \E_{x,1} (T_{2k+2}\wedge   \tau-T_{2k+1}\wedge   \tau ) ((2r_-+d)+ \epsilon)
 \nonumber \\ \nonumber \\ 
\le \mathbb E_{x,1} 1(\tau > T_{2k+1}) \underline{\lambda}_1^{-1}((2r_++d)+ \epsilon).
\end{eqnarray*}
\end{Corollary}

\noindent
{\em Proof of lemma \ref{lem2}.}
{\bf 1.} Recall that $T_0=0$ under the condition $Z_0=0$. 
We have, 
$$
T_{2k+1} = \inf(t>T_{2k}: Z_t=1).
$$
In other words, the moment $T_{2k+1}$ may be treated as ``$T_{1}$ after $T_{2k}$''. Under $Z_0=0$ the process $X_t$ coincides with $X^0_t$ until the moment $T_1$. Hence, we have on $t\in [0,T_1]$ by It\^o's formula 
\begin{eqnarray*}
dX_t^2 - 2X_t dW_t = (2X_t b_-(X_t) + d)\,dt \le (- 2r_- + d)dt,
\end{eqnarray*}
on the set $(|X_t|> M)$ due to the assumptions (\ref{b}). Further, since $1(|X_t| > M) = 1 - 1(|X_t| \le M)$, we obtain
\begin{eqnarray*}
\int_0^{T_1\wedge   \tau} 2X_t b_-(X_t)dt 
 \\\\
= 
\int_0^{T_1\wedge   \tau} 2X_t b_-(X_t) 1(|X_t| > M)dt  
+\int_0^{T_1\wedge   \tau} 2X_t b_-(X_t)1(|X_t| \le M)dt 
 \\\\
\le -  2r_- \int_0^{T_1\wedge   \tau}1(|X_t| > M)dt  
+\int_0^{T_1\wedge   \tau} 2M \|b\| 1(|X_t| \le M)dt 
 \\\\ 
= -  2r_- \int_0^{T_1\wedge   \tau}1dt  
+\int_0^{T_1\wedge   \tau} (2M \|b\| + 2r_-) 1(|X_t| \le M)dt 
 \\\\ 
\le  -  2r_- \int_0^{T_1\wedge   \tau}1dt  
+(2M \|b\| +2r_-) \int_0^{T_1\wedge   \tau} 1(|X_t| \le M)dt.
\end{eqnarray*}
Thus, always for $|x|>M_1$, 
\begin{eqnarray*}
\mathbb E_{x,z} \int_0^{T_1\wedge   \tau} 2X_t b_-(X_t)dt  
 \\\\
\le -  2r_- \mathbb E_{x,z} \int_0^{T_1\wedge   \tau}1dt
+ (2M \|b\| + 2r_-) \mathbb E_{x,z}\int_0^{T_1\wedge   \tau} 1(|X_t| \le M)dt 
\\\\
= -  2r_- \mathbb E\int_0^{T_1\wedge   \tau}1dt
+ (2M \|b\| +2r_-)\mathbb  E_{x,z}\int_0^{T_1\wedge   \tau} 1(|X_t| \le M)dt 
 \\\\
\le -  2r_- \mathbb E\int_0^{T_1\wedge   \tau}1dt 
+ (2M \|b\| + 2r_-) \mathbb E_{x,z}\int_0^{T_1} 1(|X_t| \le M)dt 
  \\\\
\le - 2 r_- \mathbb E\int_0^{T_1\wedge   \tau}1dt 
+(2M \|b\| +2 r_-) \delta .
\end{eqnarray*}
For a fixed $\epsilon>0$ let us choose $\delta  = \bar\lambda_0^{-1}\epsilon / (2M \|b\| + 2r_-) $. Then, since $|x|>M_1$ implies $T_1\wedge   \tau = T_1$  on $(Z_0=0)$, and since 
\begin{equation}\label{et1}
\bar{\lambda}_0^{-1} \le 
\mathbb E_{x,0} T_1 \le \underline{\lambda}_0^{-1}, 
\end{equation}
we get with $z=0$
\begin{eqnarray*}
\mathbb E_{x,z} X_{T_1\wedge   \tau}^2 - x^2  
\le - (2r_--d)\mathbb E_{x,z} \int_0^{T_1} dt  
+ \overline{\lambda}_0^{-1}\epsilon 
 \\\\
=- (2r_--d)\mathbb E_{x,z} {T_1}  
+ \overline{\lambda}_0^{-1}\epsilon 
\stackrel{(\ref{et1})}\le - \overline{\lambda}_0^{-1}((2r_--d)- \epsilon). 
\end{eqnarray*}
Substituting here $x$ by $X_{T_{2k}}$ and writing $\mathbb E_{x,z}(\cdot |{\mathcal F}_{T_{2k}})$ instead of $\mathbb E_{x,z}(\cdot)$, and multiplying by $1(\tau > T_{2k})$, we obtain the bounds (\ref{ele2a}) and (\ref{ele2a0}), as required. 

Note that the bound (\ref{et1}) follows straightforwardly from 
\begin{align*}
\mathbb E_{x,0} T_1 = \int_0^\infty \mathbb P_{x,0}(T_1 \ge t)dt 
= \int_0^\infty \mathbb E_{x,0} \mathbb P_{x,0}(T_1 \ge t | {\cal F}^{X^0}_{t})dt
 \\\\
= \mathbb E_{x,0}\int_0^\infty \exp(-\int_0^t \lambda_0(X^0_s)ds)dt
\le \int_0^\infty \exp(-\int_0^t \underline{\lambda}_0 ds)dt
 \\\\
= \int_0^\infty \exp(-t \underline{\lambda}_0)dt 
= \underline{\lambda}_0^{-1}, 
\end{align*}
and similarly
\begin{align*}
\mathbb E_{x,0} T_1 
= \int_0^\infty \mathbb E_{x,0} \exp(-\int_0^t \lambda_0(X^0_s)ds)dt
 \\\\
\ge \int_0^\infty \exp(-\int_0^t \overline{\lambda}_0 ds)dt
= \int_0^\infty \exp(-t \overline{\lambda}_0)dt = \overline{\lambda}_0^{-1}. 
\end{align*}

~

\noindent
{\bf 2.} The condition $Z_0=1$ implies the inequality $T_0>0$. 
We have, 
$$
T_{2k+2} = \inf(t>T_{2k+1}: Z_t=0).
$$
In other words, the moment $T_{2k+2}$ may be treated as ``$T_{0}$ after $T_{2k+1}$''. Under the condition $Z_0=1$ the process $X_t$ coincides with $X^1_t$ until the moment $T_0$. Hence, we have on $[0,T_0]$ by It\^o's formula 
\begin{eqnarray*}
dX_t^2 - 2X_t dW_t = 2X_t b_+(X_t)dt + dt \le (2r_+ + d)dt,
\end{eqnarray*}
on the set $(|X_t|> M)$ due to the assumptions (\ref{b}). Further, since $1(|X_t| > M) = 1 - 1(|X_t| \le M)$, we obtain
\begin{eqnarray*}
\int_0^{T_0\wedge   \tau} 2X_t b_+(X_t)dt  
 \\\\
= 
\int_0^{T_0\wedge   \tau} 2X_t b_+(X_t) 1(|X_t| > M)dt  
+\int_0^{T_0\wedge   \tau} 2X_t b_+(X_t)1(|X_t| \le M)dt 
 \\\\
\le 2r_+ \int_0^{T_0\wedge   \tau}1(|X_t| > M)dt  
+\int_0^{T_0\wedge   \tau} 2M \|b\| 1(|X_t| \le M)dt 
 \\\\ 
= 2r_+ \int_0^{T_0\wedge   \tau}1dt  
+\int_0^{T_1\wedge   \tau} (2M \|b\| - 2r_+) 1(|X_t| \le M)dt 
 \\\\ 
\le  2r_+ \int_0^{T_0\wedge   \tau}1dt  
+2M \|b\| \int_0^{T_0\wedge   \tau} 1(|X_t| \le M)dt.
\end{eqnarray*}
Thus, for $|x|>M_1$ and with $z=1$ we have, 
\begin{eqnarray*}
\E_{x,z} \int_0^{T_0\wedge   \tau} 2X_t b_+(X_t)dt  
 \\\\
\le 2r_+ \E_{x,z}\int_0^{T_0\wedge   \tau}1dt
+ 2M \|b\|  E_{x,z}\int_0^{T_0\wedge   \tau} 1(|X_t| \le M)dt 
\\\\
= 2r_+ {\mathbb E}_{x,z}\int_0^{T_0\wedge   \tau}1dt
+ 2M \|b\| {\mathbb  E}_{x,z}\int_0^{T_1\wedge   \tau} 1(|X_t| \le M)dt 
 \\\\
\le 2r_+ {\mathbb E}_{x,z}\int_0^{T_0\wedge   \tau}1dt 
+ 2M \|b\|  {\mathbb E}_{x,z}\int_0^{T_0} 1(|X_t| \le M)dt 
  \\\\
\le 2 r_+ \mathbb E_{x,z}\int_0^{T_0\wedge   \tau}1dt 
+ 2M \|b\|\delta .
\end{eqnarray*}
For a fixed $\epsilon>0$ let us choose $\delta  = \underline{\lambda}_1^{-1}\epsilon / (2M \|b\|) $. Then, since $|x|>M_1$ implies $T_0\wedge   \tau = T_0$ on the set $(Z_0=1)$, we get (recall that $z=1$)
\begin{eqnarray*}
\mathbb E_{x,z} X_{T_0\wedge   \tau}^2 - x^2  
\le  (2r_++d)\mathbb E_{x,z} \int_0^{T_0} dt  
+ \underline{\lambda}_1^{-1}\epsilon
 \\\\
= (2r_++d)\mathbb E_{x,z} {T_0}  
+ \underline{\lambda}_1^{-1}\epsilon
\le \underline{\lambda}_1^{-1}((2r_++d)+ \epsilon). 
\end{eqnarray*}
Substituting here $X_{T_{2k+1}}$ instead of $x$ and writing $\mathbb E_{x,z}(\cdot |{\mathcal F}_{T_{2k+1}})$ instead of $\mathbb E_{x,z}(\cdot)$, and multiplying by $1(\tau > T_{2k+1})$, we obtain the bounds  (\ref{ele2b}) and (\ref{ele2b0}), as required. Lemma \ref{lem2} is proved. \hfill 
{\em QED} 

~

\noindent
{\em Proof of corollary \ref{Cor2}} is straightforward by taking expectations.

\begin{lemma}\label{lem3}
If $M_1$ is large enough, then under the assumptions of the theorem for any $k=0,1,\ldots$ 
\begin{eqnarray}\label{ele3}
\mathbb E_{x,z} (X_{T_{2k+2}\wedge   \tau}^2 |Z_0=0, {\mathcal F}_{T_{2k+1}}) 
\le  \mathbb E_{x,z} (X_{T_{2k+1}\wedge   \tau}^2  |Z_0=0, {\mathcal F}_{T_{2k+1}})  
 \nonumber \\\\
+  1(\tau > T_{2k+1}) \mathbb E_{x,z} (T_{2k+2}\wedge \tau - T_{2k+1}\wedge \tau   |Z_0=0, {\mathcal F}_{T_{2k+1}}) ) 
((2r_++1) + \epsilon))
 \nonumber \\  \nonumber\\ \label{ele30} 
\le \mathbb E_{x,z} (X_{T_{2k+1}\wedge   \tau}^2  |Z_0=0, {\mathcal F}_{T_{2k+1}})   +  1(\tau > T_{2k+1}) \underline{\lambda}_1^{-1}((2r_++1) + \epsilon)), 
\end{eqnarray}
and 
\begin{eqnarray}\label{ele7}
\mathbb E_{x,z} (X_{T_{2k+1}\wedge   \tau}^2 |Z_0=1, {\mathcal F}_{T_{2k}}) 
\le  \mathbb E_{x,z} (X_{T_{2k}\wedge   \tau}^2  |Z_0=1, {\mathcal F}_{T_{2k}})  
  \nonumber \\\\ \nonumber 
+  1(\tau > T_{2k}) \mathbb E_{x,z} (T_{2k+1}\wedge \tau - T_{2k}\wedge \tau   |Z_0=0, {\mathcal F}_{T_{2k}}) ) 
  \nonumber \\ \nonumber \\ \label{ele70}
\le  \mathbb E_{x,z} (X_{T_{2k}\wedge   \tau}^2  |Z_0=1, {\mathcal F}_{T_{2k}})  -  1(\tau > T_{2k})\overline{\lambda}_0^{-1}((2r_--1) - \epsilon)). 
\end{eqnarray}
\end{lemma}

\begin{Corollary}\label{Cor3}
If $M_1$ is large enough, then under the assumptions of the theorem for any $k=0,1,\ldots$ 
\begin{eqnarray*}
\mathbb E_{x,0} X_{T_{2k+2}\wedge   \tau}^2 - \mathbb E_{x,0} X_{T_{2k+1}\wedge   \tau}^2  
 \nonumber \\\\
\le  \E_{x,0} 1(\tau > T_{2k+1}) (T_{2k+2}\wedge \tau - T_{2k+1}\wedge \tau) 
((2r_++1) + \epsilon)
 \nonumber \\\\
=  \E_{x,0} (T_{2k+2}\wedge \tau - T_{2k+1}\wedge \tau) 
((2r_++1) + \epsilon)
 \nonumber \\  \nonumber\\ 
\le \mathbb E_{x,0}  1(\tau > T_{2k+1}) \underline{\lambda}_1^{-1}((2r_++1) + \epsilon)), 
\end{eqnarray*}
and 
\begin{eqnarray*}
\mathbb E_{x,1} X_{T_{2k+1}\wedge   \tau}^2 - 
\mathbb E_{x,1} X_{T_{2k}\wedge   \tau}^2  
  \nonumber \\\\ \nonumber 
\le \mathbb E_{x,1} 1(\tau > T_{2k}) (T_{2k+1}\wedge \tau - T_{2k}\wedge \tau) 
  \nonumber \\ \nonumber \\ 
\le   -  \mathbb E_{x,1}1(\tau > T_{2k})\overline{\lambda}_0^{-1}((2r_--1) - \epsilon)). 
\end{eqnarray*}
\end{Corollary}

\noindent
{\em Proof of lemma \ref{lem3}.} 
Let  $Z_0=0$; recall that it implies $T_0=0$. 
If $\tau \le T_{2k+1}$, then (\ref{ele3}) is trivial. Let $\tau > T_{2k+1}$.
We will substitute $x$ instead of $X_{T_{2k}}$ for a while, and will be using the solution $X^1_t$ of the equation
\begin{eqnarray}\label{sde1}
dX^1_{t} =b(X^1_{t}, 1)\, dt+ 
dW_{t}, \quad t\ge T_1,
\quad X^1_{T_1} =X_{T_1}. \nonumber
\end{eqnarray}
For $M_1$ large enough, since $|x|\wedge |X_{T_1}|>M_1$ implies $T_2\le \tau$,  and due to the assumptions (\ref{b})  the double bound
\begin{eqnarray*}
1(|X_{T_1}|>M_1) (\mathbb E_{X_{T_1},1} X_{T_2\wedge   \tau}^2 - X_{T_1\wedge   \tau}^2)  
 \\\\
\le 1(|X_{T_1}|>M_1)(\mathbb E_{X_{T_1},1} (T_2 - T_1)((2r_++d) + \epsilon)) 
 \nonumber \\\\ \nonumber
\le + 1(|X_{T_1}|>M_1)(\underline{\lambda}_1^{-1}((2r_++d) + \epsilon))
\end{eqnarray*}
is guaranteed 
in the same way as the bounds (\ref{ele2b}) and (\ref{ele2B0}) in the previous lemma. 
In particular, it follows that for $|x|>M_1$
\begin{eqnarray*}
(\mathbb E_{X_{T_1},1} X_{T_2\wedge   \tau}^2 - X_{T_1\wedge   \tau}^2)  
\le 1(|X_{T_1}|>M_1)(\mathbb E_{X_{T_1},1} (T_2\wedge\tau - T_1\wedge\tau)((2r_++d) + \epsilon))
 \nonumber \\\\ \nonumber
= + 1(|X_{T_1}|>M_1)(\underline{\lambda}_1^{-1}((2r_++d) + \epsilon)), 
\end{eqnarray*}
since $|X_{T_1}|\le M_1$ implies $\tau \le T_1$ and $\mathbb E_{X_{T_1},1} X_{T_2\wedge   \tau}^2 - X_{T_1\wedge   \tau}^2=0$.
So, on the set $|x|>M_1$ we have with $z=0$ 
\begin{eqnarray*}
\mathbb E_{x,z} (\mathbb E_{X_{T_1},1} X_{T_2\wedge   \tau}^2 - X_{T_1\wedge   \tau}^2)  
\\\\
\le 
\mathbb E_{x,z} 1(|X_{T_1}|>M_1)(\mathbb E_{X_{T_1},1} (T_2\wedge\tau - T_1\wedge\tau)
((2r_++d) + \epsilon)
 \\\\
\le 
\mathbb E_{x,z} 1(|X_{T_1}|>M_1)(\underline{\lambda}_1^{-1}((2r_++1) + \epsilon)) 
\le \underline{\lambda}_1^{-1}((2r_++d) + \epsilon).
\end{eqnarray*}
Now substituting back $X_{T_{2k}}$ in place of $x$ and multiplying by $1(\tau > T_{2k+1}) $, we obtain the inequalities (\ref{ele3}) and (\ref{ele30}), as required.

\medskip

\noindent
For $Z_0=1$ we have $T_0>0$, and the bounds (\ref{ele7}) and (\ref{ele70}) follow in a similar way. Lemma \ref{lem3} is proved.  \hfill {\em QED}

~

\noindent
{\em Proof of corollary \ref{Cor3}} is straightforward by taking expectations.

\begin{lemma}\label{lem4}
Under the assumptions of the theorem for any $k=0,1,\ldots$
\begin{align*}
 1(\tau>T_{2k+1}) \mathbb E_{X_{T_{2k+1}},1} (T_{2k+2}\wedge \tau - T_{2k+1}\wedge \tau) \ge  1(\tau>T_{2k+1})  \bar \lambda_1^{-1},
\end{align*}
and 
\begin{align*}
 1(\tau>T_{2k}) \mathbb E_{X_{T_{2k}},0} (T_{2k+1}\wedge \tau - T_{2k}\wedge \tau) \le  1(\tau>T_{2k})    \underline\lambda_0^{-1},
\end{align*}
\end{lemma}

\noindent
{\em Proof of lemma \ref{lem4}.}
On the set $\tau>T_{2k+1}$ we have, 
\begin{align*}
 \mathbb E_{X_{T_{2k+1}},1} (T_{2k+2}\wedge \tau - T_{2k+1}\wedge \tau) 
= \mathbb E_{X_{T_{2k+1}},1} (T_{2k+2} - T_{2k+1}) 
\in [\bar \lambda_1^{-1} ,\underline \lambda_1^{-1}].
\end{align*}
Similarly, on the set $\tau>T_{2k}$
\begin{align*}
 \mathbb E_{X_{T_{2k}},0} (T_{2k+1}\wedge \tau - T_{2k}\wedge \tau) 
= \mathbb E_{X_{T_{2k}},0} (T_{2k+1} - T_{2k}) 
\in [\bar \lambda_0^{-1} ,\underline \lambda_0^{-1}].
\end{align*}
On the sets $\tau\le T_{2k+1}$ and $\tau\le T_{2k}$, respectively, both sides of the required inequalities equal zero.  Lemma \ref{lem4} follows. \hfill QED

\section{Proof of theorem 1}
Consider the case $Z_0=0$ where $T_0=0$. 
Since the identity
$$
\tau\wedge T_n = \tau\wedge T_0 + \sum_{m=0}^{n-1} 
((T_{m+1}\wedge   \tau) - (T_{m}\wedge   \tau))
$$
we have, 
$$
\mathbb E_{x,z}(\tau\wedge T_n) = \mathbb E_{x,z}\tau\wedge T_0 + \mathbb E_{x,z}\sum_{m=0}^{n-1} 
((T_{m+1}\wedge   \tau) - (T_{m}\wedge   \tau)),
$$
Due to the convergence $T_n\uparrow \infty$, we get by the monotone convergence theorem
\begin{align}\label{etau}
\mathbb E_{x,z} \tau = \mathbb E_{x,z}\tau\wedge T_0 + \sum_{m=0}^{\infty} 
\mathbb E_{x,z} ((T_{m+1}\wedge   \tau) - (T_{m}\wedge   \tau))
  \\ \nonumber \\ \nonumber 
= \mathbb E_{x,z}\tau\wedge T_0 + \sum_{k=0}^{\infty} 
\mathbb E_{x,z} ((T_{2k+1}\wedge   \tau) - (T_{2k}\wedge   \tau)) 
 \\\nonumber \\\nonumber 
+ \sum_{k=0}^{\infty} 
\mathbb E_{x,z} ((T_{2k+2}\wedge   \tau) - (T_{2k+1}\wedge   \tau)).
\end{align}
By virtue of the corollary \ref{Cor2}, 
we have
\begin{align*}
\mathbb E_{x,z} (T_{2k+1}\wedge \tau - T_{2k}\wedge \tau) 
\le ((2r_--d)- \epsilon)^{-1}\left(\mathbb E_{x,z} X_{T_{2k+1}\wedge   \tau}^2 - \mathbb E_{x,z} X_{T_{2k}\wedge   \tau}^2\right). 
\end{align*}
Therefore, 
\begin{align*}
\mathbb E_{x,0} X_{T_{2m+2}\wedge   \tau}^2 - x^2 
 \\\\
\le ((2r_++d)+ \epsilon) \sum_{k=0}^{m} \mathbb E_{x,0} (T_{2k+2}\wedge \tau - T_{2k+1}\wedge \tau)
 \\\\
- ((2r_--d)- \epsilon) \sum_{k=0}^{m} \mathbb E_{x,0} (T_{2k+1}\wedge \tau - T_{2k}\wedge \tau) 
 \\\\
= \sum_{k=0}^{m} \left(- ((2r_--d)- \epsilon) (\mathbb E_{x,0} (T_{2k+1}\wedge \tau - T_{2k}\wedge \tau) 
 \right.\\\\ \left.
+((2r_++d)+ \epsilon) \mathbb E_{x,0} (T_{2k+2}\wedge \tau - T_{2k+1}\wedge \tau) \right).
\end{align*}
By virtue of Fatou's lemma we get
\begin{align}\label{est1}
x^2 \ge ((2r_--d)- \epsilon)  \sum_{k=0}^{m}  (\mathbb E_{x,0} (T_{2k+1}\wedge \tau - T_{2k}\wedge \tau) 
 \nonumber \\\\ \nonumber
- ((2r_++d)+ \epsilon)   \sum_{k=0}^{m} \mathbb E_{x,0} (T_{2k+2}\wedge \tau - T_{2k+1}\wedge \tau).
\end{align}
Note  that $1(\tau > T_{2k+1}) \le 1(\tau > T_{2k})$. 
So, $\mathbb P_{x,0}(\tau > T_{2k+1}) \le \mathbb P_{x,0}(\tau > T_{2k})$. Hence, 
\begin{align*}
\overline\lambda_0 \mathbb E_{x,0} (T_{2k+1}\wedge \tau - T_{2k}\wedge \tau)
- \underline\lambda_1 \mathbb E_{x,0} (T_{2k+2}\wedge \tau - T_{2k+1}\wedge \tau) 
 \\\\
= \overline\lambda_0 \mathbb E_{x,0} (T_{2k+1}\wedge \tau - T_{2k}\wedge \tau)1(\tau\ge T_{2k})
 \\\\
- \underline\lambda_1 \mathbb E_{x,0} (T_{2k+2}\wedge \tau - T_{2k+1}\wedge \tau) 1(\tau\ge T_{2k+1})
 \\\\
= \overline\lambda_0 \mathbb E_{x,0} 1(\tau > T_{2k}) \mathbb E_{X_{T_{2k}}}(T_{2k+1}\wedge \tau - T_{2k}\wedge \tau)
 \\\\
- \underline\lambda_1 \mathbb E_{x,0} 1(\tau > T_{2k+1}) \mathbb E_{X_{T_{2k+1}}}(T_{2k+2}\wedge \tau - T_{2k+1}\wedge \tau) 
 \\\\
\ge \overline{\lambda}_0 \mathbb E_{x,0} 1(\tau > T_{2k}) \overline\lambda_0^{-1}
- \underline\lambda_1 \mathbb E_{x,0} 1(\tau > T_{2k+1}) \underline\lambda_1^{-1}
 \\\\
= \mathbb E_{x,0} 1(\tau > T_{2k}) 
- \mathbb E_{x,0} 1(\tau > T_{2k+1}) \ge 0. 
\end{align*}
Thus, 
$$
\mathbb E_{x,0} (T_{2k+2}\wedge \tau - T_{2k+1}\wedge \tau)  
\le \frac{\overline\lambda_0}{\underline\lambda_1} \mathbb E_{x,0} (T_{2k+1}\wedge \tau - T_{2k}\wedge \tau).
$$
Therefore, we estimate
\begin{align*}
((2r_++d)+ \epsilon)   \sum_{k=0}^{m} \mathbb E_{x,0} (T_{2k+2}\wedge \tau - T_{2k+1}\wedge \tau) 
 \\\\
\le ((2r_++d)+ \epsilon) \frac{\overline\lambda_0}{\underline\lambda_1} \sum_{k=0}^{m} \mathbb E_{x,0} (T_{2k+1}\wedge \tau - T_{2k}\wedge \tau) 
 \\\\
= q  ((2r_--d)- \epsilon) \sum_{k=0}^{m} \mathbb E_{x,0} (T_{2k+1}\wedge \tau - T_{2k}\wedge \tau).
\end{align*}
So, (\ref{est1}) implies that 
\begin{align*}\label{est2}
x^2 \ge ((2r_--d)- \epsilon)  \sum_{k=0}^{m}  (\mathbb E_{x,0} (T_{2k+1}\wedge \tau - T_{2k}\wedge \tau) 
 \nonumber \\\\ \nonumber
- ((2r_++d)+ \epsilon)   \sum_{k=0}^{m} \mathbb E_{x,0} (T_{2k+2}\wedge \tau - T_{2k+1}\wedge \tau) 
 \\\\
\ge  (1-q) ((2r_--d)- \epsilon)  \sum_{k=0}^{m}  (\mathbb E_{x,0} (T_{2k+1}\wedge \tau - T_{2k}\wedge \tau) 
 \nonumber \\\\ \nonumber
 \ge \frac{1-q}2\, ((2r_--d)- \epsilon)  \sum_{k=0}^{m}  (\mathbb E_{x,0} (T_{2k+1}\wedge \tau - T_{2k}\wedge \tau) 
 \\\\
+ \frac{1-q}{2q}\, ((2r_++d)+ \epsilon)  \sum_{k=0}^{m} \mathbb E_{x,0} (T_{2k+2}\wedge \tau - T_{2k+1}\wedge \tau).
\end{align*}
Denoting $\displaystyle c:= \min\left(\frac{1-q}{2q}\, ((2r_++d)+ \epsilon), \frac{1-q}2\, ((2r_--d)- \epsilon)\right)$, we conclude that 
\begin{align*}
x^2 \ge c \sum_{k=0}^{2m}  \mathbb E_{x,0} (T_{k+1}\wedge \tau - T_{k}\wedge \tau).
\end{align*}
So, as $m\uparrow \infty$, by the monotone convergence theorem we get the inequality 
\begin{align*}
\sum_{k=0}^{\infty}  \mathbb E_{x,0} (T_{k+1}\wedge \tau - T_{k}\wedge \tau) \le c^{-1}x^2.
\end{align*}
Due to (\ref{etau}), it implies that 
\begin {equation}\label{simpleest}
\mathbb E_{x,0} \tau \le c^{-1}x^2,
\end{equation}
as required. Recall that this bound is established for $|x|>M_1$, while in the case of $|x|\le M_1$ the left hand side in this inequality is just zero. 

~

\noindent
In the case of $Z_0=1$ (and, hence, $T_0>0$), we have to add the value $\mathbb E_{x,z} T_0$ satisfying the bound $\mathbb E_{x,1} T_0\le \underline\lambda_1^{-1}$ to the right hand side of (\ref{simpleest}), which leads to the bound  (\ref{e3}), as required. Theorem 1 is proved. \hfill {\em QED}

~

\noindent
\begin{Remark}
In turn, positive recurrence for the model under the consideration implies existence of the invariant measure, see  
\cite[Section 4.4]{Khasminskii}. Convergence to this invariant measure in total variation and, hence, uniqueness of this measure follows, for example,  due to the  coupling method in a standard way. The full proof of this corollary will be presented in the publications to follow.
\end{Remark}

~

\noindent
\begin{Remark}
The results of the paper may be extended to the equation 
\begin{align}\label{sdes2}
dX_{t} =b(X_{t}, Z_t)\, dt+ \sigma (X_{t}, Z_t)\,dW_{t}, \quad t\ge 0, 
\quad X_{0} =x, \; Z_0=z,
\end{align}
with a Borel measurable $\sigma$ 
under the assumptions of the existence of a strong solutions, or of a weak solution which is weakly unique (because the strong Markov property is needed), in addition to the standing  balance type conditions replacing (\ref{b}) and (\ref{c1}) (while (\ref{al}) is still valid): $a(x,z)=\sigma\sigma^*(x,z)$ and

\begin{equation}\label{b2}
2x b(x,0) + \mbox{Tr}\,(a(x,0)) \le - R_-, \;\; 2x b(x,1) + \mbox{Tr}\,(a(x,1))\le + R_+, \;\; \forall \,|x|\ge M, 
\end{equation}
with some $R_-, R_+ >0$, and 
\begin{equation}\label{c12}
\underline{\lambda}_1 R_- > \overline{\lambda}_0 R_+,
\end{equation}
where the definitions of $\underline{\lambda}_1$ and $\overline{\lambda}_0 $ do not change. 
The proofs now involve the diffusion coefficient and use the assumptions (\ref{b2}) and (\ref{c12}), but otherwise remain the same as in the case of the unit diffusion matrix.
\end{Remark}


~

\noindent
{\bf Availability of data and material.} Data sharing is not applicable to this article as no datasets were generated or analysed during the current study.

\end{document}